\documentclass[runningheads]{llncs}
\usepackage{framed}
\usepackage[T1]{fontenc}
\usepackage{graphicx}
\usepackage{hyperref}
\usepackage{color}

\usepackage{amsmath,amssymb,amsfonts, enumitem, fancyhdr, color, comment, graphicx, environ, tikz-cd, physics}

\newcommand{\mb}[1]{\mathbf{#1}}
\newcommand{\mc}[1]{\mathcal{#1}}

\newcommand{\cx}{\mathbb{C}}
\newcommand{\re}{\mathbb{R}}

\newcommand{\proj}{\mathbb{P}}

\newcommand{\im}{\operatorname{im}}
\newcommand{\p}{\partial}

\newcommand\setb[1]{\left\{#1\right\}}
\newcommand{\vbrack}[1]{\left\langle #1 \right\rangle}
\renewcommand{\bar}[1]{\overline{#1}}

\renewcommand\vec[1]{\mb #1}

\usepackage{listings}
\definecolor{codedarkgreen}{RGB}{51, 133, 4}
\definecolor{codemaroon}{RGB}{133, 5, 63}
\definecolor{codeteal}{RGB}{0, 145, 109}
\definecolor{codepurple}{RGB}{123, 35, 125}
\lstdefinelanguage{Macaulay2}{
    basicstyle=\normalsize\ttfamily,
    alsoletter=",
    classoffset=1,
    keywords={coefficients,toList,transpose,det,factor,netList,subsets,genericMatrix,needsPackage,presentation,generators,gens,selectInSubring,i1,i2,i3,i4,i5,i6,i7,i8,i9,i10,i11,i12,i13,i14,i15,i16,i17,i18,i19,i20,i21,i22,i23,i24,i25,flatten,ideal,entries,basis,coefficient,apply,kernel,map,minors,random,drop},
    keywordstyle={\color{blue}},
    classoffset=2,
    breaklines=true,
    morekeywords={"SubalgebraBases","Engine","Brackets"},
    keywordstyle={\color{codemaroon}},
    classoffset=3,
    morekeywords={QQ,CacheTable,Matrix,Eliminate},
    keywordstyle={\color{codedarkgreen}},
    classoffset=4,
    morekeywords={restart,false,true,Weights,Limit,Lex,MonomialOrder,CoefficientRing},
    keywordstyle={\color{codeteal}},
    classoffset=5,
    morekeywords={list,for,in,from,to,of},
    keywordstyle={\color{codepurple}},
    xleftmargin=1em,
    xrightmargin=1em,
    columns=fullflexible,
    keepspaces=true,
    stepnumber=1,
    numbers=none,
    captionpos=b,
    showspaces=false,
    frame=none
}

\newcommand{\Sym}{\operatorname{Sym}}
\newcommand{\sing}{\operatorname{sing}}
\newcommand{\augmentedJacobian}{\widehat{\mc{J}}}

\providecommand{\demph}[1]{\emph{#1}}

\begin{document}
\title{Euclidean Distance Degrees in Macaulay2}
%
%
\author{William Huang\inst{1}\orcidID{0000-0001-6287-772X} \and
Jose Israel Rodriguez\inst{2}\orcidID{0000-0003-3140-9944} 
}
\authorrunning{W. Huang and J. Rodriguez}
%
\institute{University of British Columbia, Vancouver, Canada 
\and
University of Wisconsin --- Madison, Madison, Wisconsin\\
\email{jose@math.wisc.edu}}
\maketitle              
\begin{abstract}
We introduce \texttt{EuclideanDistanceDegree}, a \emph{Macaulay2} package that implements symbolic and numerical methods for computing Euclidean Distance (ED) degrees. The package includes symbolic methods based on minors and conormal varieties as well as numerical methods for unit and generic ED degrees using tools from numerical algebraic geometry. We illustrate the package functionality with a range of examples in the paper and in the accompanying GitHub repository.

\end{abstract}

\section{Introduction}\label{sec:introduction}

One of the most widely studied invariants in algebraic optimization is the Euclidean Distance (ED) degree. It measures the complexity of nearest-point problems for real algebraic varieties: for a generic data point, it is the number of critical points of the squared Euclidean distance function on the variety. As many models across science and engineering are expressed in terms of polynomial systems, tools for computing ED degrees of arbitrary varieties allow researchers to better understand the algebraic complexity of their models.

In the \texttt{EuclideanDistanceDegree} package for the computer algebra system \emph{Macaulay2} \cite{M2}, we implement methods for computing these degrees. These include symbolic methods to compute exact ED degrees for affine models as well as numerical methods using homotopy continuation. This package relies on the Bertini software package \cite{BHSW06} and the corresponding \emph{Macaulay2} wrapper \cite{bates2013bertinimacaulay2}. We also utilize the \texttt{MonodromySolver} package \cite{MonodromySolverSource} to compute ED degrees of parameterized varieties.

\section{Computing ED Degrees}\label{sec:background}

The Euclidean Distance degree arises from the following problem: given a variety $X \subseteq \re^n$ and data point $\vec u \in \re^n \setminus X$, find the point $\vec x^* \in X$ which minimizes the squared distance
\begin{equation}
    \|\vec x - \vec u \|^2 = \sum_{i= 1}^n (x_i - u_i)^2
\end{equation}
We also consider the weighted squared distance, which depends on a weight vector $\vec w \in \re^n_{>0}$ 
in addition to the data point:
\begin{equation}
    \|\vec x - \vec u\|^2_{\vec w} = \sum_{i= 1}^n w_i(x_i - u_i)^2.
\end{equation}
By viewing $X$ as a complex variety in $\cx^n$, the number of nonsingular complex critical points of $ \|\vec x - \vec u\|^2_{\vec w}$ is the \demph{weighted Euclidean distance degree} for the variety $X$ with respect to the weights $\vec w$. This number is independent of the choice of data point $\vec u$, provided it is sufficiently general \cite{Draisma2016}. We refer to the choices of arbitrary weights, unit weights, and generic weights as the weighted, unit, and generic ED degrees, respectively.

\subsection{Minors Method}\label{subsec:symbolic}
Given an irreducible real algebraic variety $X$, the ED degree can be directly computed by viewing $X$ as a complex variety. Let $I_X = \vbrack{f_1, \ldots, f_k}$ be its prime ideal and $\mc J = (\p f_i/\p x_j)$ the $k \times n$ Jacobian matrix. Let $c$ be the codimension of $X$. Then the singular locus $X_{\sing}$ of $X$ is defined by the ideal
\begin{equation}
    I_{X_{\sing}} = I_X + \vbrack{\text{$c \times c$ minors of $\mc J$}}
\end{equation}
The augmented Jacobian matrix $\augmentedJacobian$ is the $(k+1) \times n$ matrix obtained by placing the row vector $(w_1(x_1 - u_1), \ldots, w_n(x_n - u_n))$ above the Jacobian, i.e.
\begin{equation}
    \augmentedJacobian = \mqty(\vec w \cdot (\vec x - \vec u) \\ \mc J).
\end{equation}
The weighted critical ideal of $X$ with respect to the data point $\vec u$ and weights $\vec w$ is the following saturation:
\begin{equation}\label{eq:critideal}
    \mc C_{X, \vec u, \vec w} = \qty(I_X + \text{$(c+1)\times (c+1)$ minors of $\augmentedJacobian$}) : \qty(I_{X_{\sing}})^\infty
\end{equation}
For general data points $\vec u$, the variety defined by this ideal is finite and it consists of the nonsingular critical points \cite[Lemma 2.1]{Draisma2016} of the squared weighted distance function. Thus the degree of $\mc C_{X, \vec u, \vec w}$ is exactly the weighted ED degree of $X$ with respect to the weights $\vec w$. This method of computing ED degrees can be implemented symbolically and we call it the \demph{minors method}.

\subsection{Left Kernel Method}\label{subsec:leftkernel}
For varieties of large codimension, the minors method is computationally intensive due to the saturation step along with the possibility of there being exponentially many minors. To avoid incurring these costs, we introduce the \demph{left kernel method}, which uses Lagrange multipliers to solve the equations defining critical points. We present this method for unit ED degrees, but the same ideas will apply to the weighted case.

Let $X$ be an irreducible real algebraic variety defined by the polynomials $f_1, \ldots, f_k$ and denote $f(\vec x) = \mqty(f_1(\vec x) & \cdots & f_k(\vec x))^T$. By introducing Lagrange multipliers $\lambda = [\lambda_0: \lambda_1: \cdots : \lambda_k] \in \proj^k$, we define
\begin{equation}\label{eq:G}
    G(\vec x, \lambda) = \begin{pmatrix}
        f(\vec x) \\
        \lambda_0(\vec x - \vec u) + \sum_{i=1}^k \lambda_i \nabla f_i(\vec x)^T
    \end{pmatrix}
\end{equation}

The system $G(\vec x, \lambda) = 0$ consists of $n+k$ equations in $n+k$ degrees of freedom. If we assume that $X$ is a complete intersection given by the polynomials $f_1, \ldots, f_k$, then for every pair $(\vec x, \lambda) \in V(G(\vec x, \lambda))$, the vector $\vec x$ is a critical point of the squared distance function \cite{Hauenstein2013}. Thus the number of isolated solutions to $G(\vec x, \lambda)$ is exactly the unit ED degree of $X$. In our implementation, we use the Bertini software package to solve this system numerically.

\subsection{Homotopy Method}\label{subsec:numeric}
The left kernel method can be extended into what we call the \demph{homotopy method}, which is especially useful for exploring the behavior of weighted ED degrees for special choices of weights or data. Let $F = \setb{f_1, \ldots, f_k}$ and $ G = \setb{g_1, \ldots, g_m}$ be two sets of polynomials such that $V(G)$ is a complete intersection containing $V(F)$ as an irreducible component. Denote by $\mc J_G$ the Jacobian of $G$ and introduce Lagrange multipliers $\lambda = [\lambda_0 : \lambda_1 : \cdots : \lambda_m] \in \proj^{m}$. We seek to solve the system
\begin{equation}\label{eq:symbSystem}
    M(\vec x, \lambda) = \begin{pmatrix}
    \lambda_0 & \cdots & \lambda_m
    \end{pmatrix} \begin{pmatrix}
    \nabla_{\vec x} \| \vec x - \vec u \|_{\vec w} ^2 \\
    \mc J_G 
\end{pmatrix}
\end{equation}
subject to $\vec x \in V(F)$. The number of such solutions $(\lambda, \vec x)$ is the ED degree of $V(F)$. Using homotopy continuation, we do this computation in two stages:
\begin{enumerate}
    \item We construct a start system of the form \eqref{eq:symbSystem} using generic weights. We then leverage Bertini to solve the start system to obtain an initial bound on the ED degree of $V(F)$.
    
    \item Using the results of stage one, we perform a parameter homotopy to track the solutions from the start system to the target system. By filtering out points that correspond to infinity or do not vanish on $F$, we get the desired ED degree.
\end{enumerate}

The two-stage approach allows one to efficiently compute ED degrees for multiple choices of weights $\vec w$ or data $\vec u$ because the first stage only needs to be executed once. Correctness of the left-kernel and homotopy formulations assumes that the relevant input polynomials define a complete intersection: $F$ for the left-kernel method, and $G$ for the homotopy method.

\subsection{Conormal Varieties}\label{subsec:conormal}
Suppose instead that the variety $X$ is given by a parameterization. More precisely, let $\phi_1, \ldots, \phi_n$ be polynomials in $d$ variables. We assume $n > d$ and define the map
\begin{equation}\label{eq:param}
    \phi: \cx^d \to \cx^n \qquad \vec x \mapsto (\phi_1(\vec x), \ldots, \phi_n(\vec x))
\end{equation}
The variety $X$ is the Zariski closure of $\im \phi$. Suppose $\dim X = d$ so that the Jacobian map
\begin{equation}
    \nabla_{\vec x}(\phi_1, \ldots, \phi_n): \cx^d \to \cx^{n \times d} \qquad \vec y \mapsto \mqty(\dfrac{\p \phi_i}{\p x_j}(\vec y)) = \mc J(\vec y)
\end{equation}
is generically full rank. The points for which the Jacobian is not full rank form the set of critical points for the parameterization \eqref{eq:param}. If $\vec y \in \cx^d$ is not a critical point, then the columns of $\mc J(\vec y)$ span a $d$-dimensional subspace of $\cx^n$. We globally describe the orthogonal complement of this subspace by finding a spanning set for $\ker \mc J^T$. 

Since the kernel may be over-parameterized, we evaluate its generators at a general point $\vec y \in \cx^d$. We then find a subset of $n - d$ elements which span a $(n - d)$-dimensional subspace and select the corresponding generators of $\ker \mc J^T$. Let $M \in (\cx[\vec x])^{n \times (n-d)}$ be the matrix whose columns are those elements and define the parametric map
\begin{equation}\label{eq:conormal}
    \Psi_{M}: \cx^d \times \cx^{n-d} \to \cx^n \times \cx^n \qquad (\vec x,\lambda) \mapsto \left(\phi(\vec x), M(\vec x) \cdot \mqty(\lambda_1 \\ \vdots \\ \lambda_{n-d})\right)
\end{equation}
Under the given assumptions, the image of $\Psi_{M}$ is $n$-dimensional and we call its Zariski closure an \emph{affine conormal variety}. This construction can be used to compute ED degrees for varieties given via parameterizations. We wish to find points $\vec x$ such that
\begin{equation}
    \sum_{i=1}^n 2w_i(\phi_i(\vec x) - u_i)\frac{\p \phi_i}{\p x_j}(\vec x) = 0 \qquad j = 1, \ldots, d
\end{equation}
This is equivalent to the statement
\begin{equation}
    \mc J(\vec x)^T \mqty(2w_1(\phi_1(\vec x) - u_1) \\ \vdots \\ 2w_n(\phi_n(\vec x) - u_n)) = 0
\end{equation}
Using the parameterization \eqref{eq:conormal}, we see that critical points are exactly the pairs $(\vec x, \lambda)$ which satisfy
\begin{equation}
    \mqty(2w_1(\phi_1(\vec x) - u_1) \\ \vdots \\ 2w_n(\phi_n(\vec x) - u_n)) = M(\vec x) \cdot \mqty(\lambda_1 \\ \vdots \\ \lambda_{n-d})
\end{equation}
The ideal generated by these $n$ equations is the critical ideal of the parameterized variety $X$ and its degree is the weighted ED degree of $X$.

\section{Implementations}\label{sec:methods}

The \texttt{EuclideanDistanceDegree} package implements the four methods for computing ED degrees described in Section~\ref{sec:background}. Each method is implemented across three variants with differing levels of generality:
\begin{itemize}[nosep]
    \item \texttt{Weight}: user specified weights and data
    \item \texttt{Unit}: unit weights and random data
    \item \texttt{Generic}: random weights and data
\end{itemize}

All numerical methods rely on homotopy continuation and may underestimate the ED degree due to path failures, e.g. when paths diverge to infinity. The package writes the corresponding Bertini input files to a temporary directory. This directory can be changed using the \texttt{TempDirectory} option, allowing users to inspect, modify, and rerun Bertini computations. Users with Bertini experience may mitigate the possibility of path failures by modifying the configuration in these files, for example by changing tolerance or step sizes.

\subsection{Minors Method: Symbolic Computation}\label{subsec:determinantal}
To compute ED degrees symbolically using the minors method as described in Section~\ref{subsec:symbolic}, the package includes the three methods \texttt{symbolicWeightEDDegree}, \texttt{determinantalUnitEDDegree}, and \texttt{determinantalGenericEDDegree}. 
The latter two methods take as input a list of polynomials defining an irreducible variety and computes the corresponding ED degree. For additional control over the computation, the \texttt{symbolicWeightEDDegree} method also takes a data point and weight vector. 

One distinctive feature of these methods is the ability to return the equations defining the critical ideal \eqref{eq:critideal}. This is done by setting the optional argument \texttt{ReturnCriticalIdeal=>true}. This functionality is particularly important for replication with other software packages like \texttt{HomotopyContinuation.jl} \cite{HomotopyContinuation-dot-JL}. As an example, suppose we want to compute the ED degree of the Dingdong surface from Herwig Hauser's algebraic surfaces gallery \cite{HauserGallery}. We can compute the unit and generic ED degrees and recover the critical ideal using the following code:

\begin{center}
\begin{lstlisting}[language=macaulay2]
i1: R = QQ[x, y, z];
i2: F = {x^2 + y^2 + z^3 - z^2};
i3: determinantalUnitEDDegree F
o3: 5
i4: determinantalGenericEDDegree F
o4: 9
i5: (U, W) = ({1,12,2}, {1,1,2});
i6: symbolicWeightEDDegree(F, U, W, ReturnCriticalIdeal => true)
o6 = ideal (12x - y, 145y^2 - 435y*z + 288z^2 - 870y + 4548z - 3096, 48z^3 + 145y*z - 144z^2 + 290y - 1516z + 1032, 3y*z^2 - 6y*z - 36z^2 + 8y + 24z)
o6: Ideal of R
\end{lstlisting}
\end{center}

\subsection{Left Kernel Method: Complete Intersections}\label{sec:completeInter}
For varieties that are complete intersections, the left kernel method may be used as described in Section~\ref{subsec:leftkernel}. The package implements the left kernel method in the \texttt{leftKernelWeightEDDegree} method along with its unit and generic variants. For example, to compute the weighted ED degree of the Daisy surface using specific data and weights, we use the following code:
\begin{center}
\begin{lstlisting}[language=macaulay2]
i1: R = QQ[x, y, z];
i2: F = {(x^2 - y^3)^2 - (z^2 - y^2)^3};
i3: (U, W) = ({.12, .23, .25}, {.15, .331, .727});
i4: dir = temporaryFileName();
i5: leftKernelWeightEDDegree(F, U, W, TempDirectory => dir)
o5: 22
\end{lstlisting}
\end{center}

\subsection{Homotopy Method: Beyond Complete Intersections}

To compute the ED degree using the homotopy method as described in Section~\ref{subsec:numeric}, the package includes the method \texttt{numericWeightEDDegree} along with its unit and generic variants. These methods differ from those described in the previous section in that they take two lists of polynomials: a list $F$ which describes the variety of interest and another list $G$ such that the variety $V(F)$ is an irreducible component of $V(G)$. For example:
\begin{center}
\begin{lstlisting}[language=macaulay2]
i1: R = QQ[x, y];
i2: F = G = {x^2 + y^2 - 1};
i3: dir = temporaryFileName();
i4: numericUnitEDDegree(F, G, TempDirectory => dir)
o4: 2
\end{lstlisting}
\end{center}

These methods are all special cases of the \texttt{homotopyEDDegree} method, which is a configurable implementation of the homotopy method. 
The homotopy is done in the two stages described in Section~\ref{subsec:numeric}: 
the first stage solves the start system while the second stage verifies solutions on the target system via path tracking. Configuration options are stored in a \texttt{NumericalComputationOptions} object and keys may be modified to customize the homotopy. Currently we support performing homotopy continuation on weights and data. As an example, we present a weight homotopy computation of a unit and generic ED degree. Note how after the unit degree is computed, only the second stage needs to be executed to compute the generic ED degree:
\begin{center}
\begin{lstlisting}[language=macaulay2]
i1: R = QQ[x1,x2,x3,x4,x5,x6];
i2: F = (minors(2, genericMatrix(R,3,2)))_*;
i3: G = drop(F, -1);
i4: NCO = newNumericalComputationOptions(F, G);
i5: NCO#"TargetWeight" = apply(#gens R, i->1);
i6: homotopyEDDegree(NCO, "Weight", true, true)
o6: 2
i7: NCO#"TargetWeight" = apply(#gens R, i->random RR);
i8: homotopyEDDegree(NCO, "Weight", false, true)
o8: 10
\end{lstlisting}
\end{center}

The \texttt{NumericalComputationOptions} object allows for custom Bertini options to be specified in the \texttt{BertiniStartFiberSolveConfiguration} key. Thus path failures from the \texttt{homotopyEDDegree} method can be addressed directly in \emph{Macaulay2} rather than through the Bertini input files, unlike the methods introduced previously.

Because solutions can be tracked to different data points without recomputing the initial stage, \texttt{homotopyEDDegree} enables efficient ED degree experiments. For convenience, we include the \texttt{averageNumericEDDegree} method to streamline the process. This method performs a stage-one solve at random initial data, generates $n$ random data points using a user-provided sample generator, performs a stage-two solve for each sample, counts the number of real critical points, and returns the average. In \cite[Example 4.5]{Draisma2016}, the average real ED degree of an ellipse was computed to approach $3.05$ using $10^5$ samples. As a demonstration, we show how this experiment can be run for 100 samples using our package:
\begin{center}
\begin{lstlisting}[language=macaulay2]
i1: loadPackage "Probability";
i2: R = QQ[x,y];
i3: F = G = {x^2 + 4*y^2 - 4};
i4: Z = normalDistribution();
i5: sampleGen = () -> apply(#gens R, i -> random Z)
i6: averageNumericEDDegree(F, G, sampleGen, 100);
o6 = 3.41
o6: RR (of precision 53)
\end{lstlisting}
\end{center}

\subsection{Conormal Method: Parameterizations}

For parameterized varieties, the conormal-variety approach from Section~\ref{subsec:conormal} is implemented in the method \texttt{parametrizedWeightEDDegree} along with its unit and generic variants. These methods take as input a list of polynomials $F$ which parameterize a variety $X$ and then returns the ED degree of $X$. For example: 
\begin{center}
\begin{lstlisting}[language=macaulay2]
i1: R = QQ[x1, x2, x3, x4, x5];
i2: F = {x1^2+x4^2, x2^2+x5^2, x3^2+1, x1*x2+x4*x5, x1*x3+x4, x2*x3+x5};
i3: parameterizedGenericEDDegree F
o3: 52
\end{lstlisting}
\end{center}
By default this is done symbolically. Setting the option \texttt{UseMonodromy=>true} will instead leverage the \texttt{MonodromySolver} \cite{MonodromySolverSource} package to solve the critical equations.

\section{Examples}\label{sec:examples}

We conclude with two application areas, polynomial neural networks and multiview varieties, that illustrate the package’s utility for parameterized models arising from applied algebraic geometry.

\subsection{Polynomial Neural Networks}

As a first application, we compute ED degrees for shallow polynomial neural networks (PNNs). This follows the work done in \cite{Kubjas_2024}, which we shall briefly describe here. An $L$ layer PNN architecture is given by a vector $\vec d = (d_0, \ldots, d_L)$ which describes layer widths and a positive integer $r$ specifying
the activation degree. The PNN corresponding to the parameters $(\vec d, r)$ is the map
\begin{equation}
    p_{(\vec d, r)} = W_L \circ \sigma_{L-1} \circ \cdots \circ \sigma_1 \circ W_1: \re^{d_0} \to \re^{d_L}
\end{equation}
where $W_i \in \re^{d_i \times d_{i-1}}$ are linear maps and the activation functions $\sigma_i$ are applied coordinate-wise, i.e.
\[ \sigma_i(\vec x) = (x_1^r, \ldots, x_n^r) \]

The resulting map $p_{(\vec d, r)}$ can be thought of as an element of $\Sym_{r^{L-1}}(\re^{d_0})^{d_L}$, the space of tuples of degree $r^{L-1}$ homogeneous polynomials in $d_0$ variables. The closure of the image of the parameterization map in this space is the \demph{neurovariety} corresponding to the architecture $\vec d$. Using our package, we compute that the $\vec d = (3,1,1)$ network with square activation has generic ED degree 13. The architecture $\vec d = (3,2,1)$ also has generic ED degree 13 for square activations.

\subsection{Multiview Varieties}

Multiview varieties arise in the study of computer vision through algebraic geometry. In this setting, the process of taking a picture can be modeled as a rational map $\proj^3 \dashrightarrow \proj^2$. Here we follow the definition put forth in \cite{duff2024metricmultiviewgeometry,finkel-rodriguez-multiview}: let $\vec C = (C_1, \ldots C_n)$ be a tuple of full rank $(h + 1) \times (N + 1)$ matrices. Define the map
\begin{equation}
    \Phi_{\vec C}: \proj^N \dashrightarrow (\proj^h)^n \qquad X \mapsto (C_1X, \ldots, C_nX)
\end{equation}

We call the tuple $\vec C$ a \emph{camera arrangement}. Indeed for $N = 3$ and $h = 2$, $\Phi_{\vec C}$ is just the process of taking $n$ pictures with $n$ cameras. The \emph{(point) multiview variety} of $\proj^N$ with respect to $\vec C$ is the Zariski closure of the map $\Phi_{\vec C}$,
denoted
$\vec C \Box \proj^N := \bar{\im \Phi_{\vec C}} \subset (\proj^h)^n$.

More generally, for a subvariety $Y \subseteq \proj^N$, the \emph{(point) multiview variety anchored at $Y$} with respect to $\vec C$ is the closure
\[ \vec C \Box Y := \bar{\Phi_{\vec C}(Y)}.\]
The restriction of a multiview variety to the standard affine chart is denoted $(\vec C \Box Y )_\text{Aff}$ and is called the \emph{affine multiview variety} of $\vec C \Box Y$. We are interested in the case where $Y$ is a degree $E$ rational surface parameterized by
\begin{equation}
    f: \proj^2 \to \proj^N \qquad [s:t:u] \mapsto [f_0(s,t,u) : \cdots : f_N(s,t,u)].
\end{equation}

Using our package, we compute for $E = 1, N = 2, h = 2,3,4$ that the ED degree of the affine multiview variety for a rational surface is 8 and when $E = 2, N = 2, h = 2$, the ED degree is 12. The case for a rational curve was worked out in \cite[Theorem 2.3]{finkel-rodriguez-multiview}
where it is proved that the ED degree is $3En-2$, notably independent of $h$. Our computations here provide evidence that a similar result may hold in the case of a rational surface.

\subsection{Comparisons with Existing Software}

Several software packages are available for computing, or helping to compute, ED degrees of algebraic varieties symbolically. For instance the \emph{Macaulay2} package \texttt{AlgebraicOptimization} \cite{Harkonen2020} is a general-purpose symbolic package for algebraic optimization problems with dedicated methods for ED degree computations. These methods involve methods based on multidegrees, projections, as well as Fritz John formulations, the latter being similar to the equations we use in the left kernel formulation. The \emph{Macaulay2} package \texttt{ToricInvariants} \cite{Helmer_2018} computes the generic ED degree of a projective toric variety along with other invariants of the dual variety. It provides a specialized tool for ED degree computations in the toric setting. 
        
To our knowledge, there is no dedicated package for numerically computing ED degrees. However,  homotopy continuation solvers such as Bertini \cite{BHSW06}, PHCpack \cite{verschelde1999phcpack}, the \emph{Macaulay2} package \texttt{MonodromySolver} \cite{DHJLKS-monodromy-journal,MonodromySolverSource}, and the Julia package HomotopyContinuation.jl \cite{HomotopyContinuation-dot-JL} can compute an ED degree by forming the critical equations for the squared distance function to a generic point and counting the isolated solutions numerically. The option \texttt{ReturnCriticalIdeal=>true} (see Section~\ref{subsec:determinantal}) allows our package to interface with other software by  returning the generating polynomials of the critical ideal \eqref{eq:critideal}.

\begin{credits}

\subsubsection{Code Availability}
Code for the provided examples is available in the \texttt{Examples} folder of the GitHub repository \href{https://github.com/JoseMath/EuclideanDistanceDegree}{https://github.com/JoseMath/EuclideanDistanceDegree}.

\subsubsection{\ackname}
This research was partially supported by the Alfred P. Sloan Foundation and by the 
National Science Foundation grant DMS-2510307. We are grateful to the anonymous reviewers for their valuable comments and suggestions, which helped improve the quality of this paper.

\subsubsection{\discintname}
The authors have no competing interests to declare that are relevant to the content of this article.

\end{credits}

%
%
%
\bibliographystyle{splncs04}
\bibliography{refs}

\end{document}